\documentclass{article}

\usepackage{amsmath,amssymb,amsthm}

\usepackage[utf8]{inputenc}

\usepackage{hyperref}

\newtheorem{theorem}{Theorem}[subsection]
\newtheorem{lemma}[theorem]{Lemma}

\newtheorem{conjecture}[theorem]{Conjecture}
\theoremstyle{definition}
\newtheorem{definition}[theorem]{Definition}
\newtheorem{example}[theorem]{Example}

\theoremstyle{remark}
\newtheorem{remark}[theorem]{Remark}

\newcommand{\PP}{\mathbb P}
\newcommand{\CC}{\mathbb C}
\newcommand{\ZZ}{\mathbb Z}
\newcommand{\cO}{\mathcal O}
\newcommand{\Mbar}{\overline M}

\DeclareMathOperator{\Aut}{Aut}
\DeclareMathOperator{\Id}{Id}

\DeclareMathOperator{\ev}{ev}

\begin{document}

\title{Frobenius Manifolds near the Discriminant and Relations in the
  Tautological Ring}
\author{Felix Janda}
\maketitle
\begin{abstract}
  We give a criterion for extending a generically semisimple (not
  necessarily conformal) Frobenius manifold locally near a smooth
  point of the discriminant to a cohomological field theory.

  As an application, we show that a large set of tautological
  relations related to the Givental--Teleman classification for any
  generically semisimple cohomological field theories follow from
  Pixton's generalized Faber--Zagier relations.
\end{abstract}
\renewcommand{\thefootnote}{\fnsymbol{footnote}}
\renewcommand{\thefootnote}{\arabic{footnote}}

\section{Introduction}

\subsection{Frobenius manifolds}

Introduced by Dubrovin in \cite{Du94}, Frobenius manifolds are
particular flat Riemannian manifolds with an algebra structure on the
space of vector fields.
The most important examples come from genus zero Gromov--Witten theory
where the algebra structure is given by quantum multiplication, and
singularity theory where instead the Milnor ring with its multiplicative
structure is considered.

A point $p$ in a Frobenius manifold $M$ is called semisimple if there
exists a basis of idempotent local vector fields $\epsilon_i$ near $p$
which are orthogonal with respect to the metric $\eta$, and can be
integrated to give the so-called canonical coordinates.
The local behavior near a semisimple point is uniquely determined by
the norms $\eta(\epsilon_i, \epsilon_i)$ of the idempotents with
respect to the metric.
Dubrovin has shown that the $\eta(\epsilon_i, \epsilon_i)$ satisfy the
system of Darboux--Egoroff equations and conversely constructs from any
generic solutions to the Darboux-Egoroff system a Frobenius manifold.

It is a natural question to ask how this local classification extends
to non-semisimple points.
A first step in this direction has been taken by Hertling \cite{He99P}
who shows that the germ of an $N$-dimensional generically semisimple
Frobenius manifold (or more generally a massive $F$-manifold) near a
smooth point $p$ on the discriminant locus of a generically semisimple
Frobenius manifold is of the form $I_2(m) \times A_1^{N - 2}$ for some
integer parameter $m \ge 3$.
This means in particular that all but two idempotent vector fields
extend to $p$, and the parameter $m$ describes the behavior of the
remaining two idempotents behave near $p$.
In this article, Frobenius manifolds such that $m = 3$ will be of most
importance.

\subsection{Cohomological field theories}

Cohomological fields theories (CohFTs) as first introduced in
\cite{KoMa97} can be viewed as higher genus analogs of Frobenius
manifolds.
A CohFT $\Omega$ defined on a vector space $V$ is a system of
multilinear maps
\begin{equation*}
  \Omega_{g, n}\colon (V^*)^{\otimes n} \to H^*(\Mbar_{g, n})
\end{equation*}
to the cohomology of the moduli space of curves $H^*(\Mbar_{g, n})$
satisfying some compatibility axioms.
By restricting to genus zero, the structure of a formal neighborhood
of a point in a Frobenius manifold with tangent space $V$ can be
extracted.
Under convergence assumptions, $V$ is the tangent space of a point of
an actual Frobenius manifold and using a shifting construction
$\Omega_{g, n}$ can be extended to a family of CohFTs depending on a
point on a Frobenius manifold.
We call such a family of CohFTs a \emph{convergent CohFT}.

For a convergent CohFT $\Omega_{g, n}$, the Givental--Teleman
reconstruction \cite{Te12} of semisimple CohFTs gives a way to
reconstruct the CohFT $\Omega_{g, n}$ from its underlying Frobenius
manifold $M$ near a semisimple point of $M$ up to a choice of
integration constants.
In many applications, $\Omega_{g, n}$ is in addition homogeneous with
respect to an Euler vector field, in which case, there is a unique
choice for the integration constants.

The Givental--Teleman reconstruction gives a way to extend the
neighborhood of a semisimple point of a Frobenius manifold to a
convergent CohFT.
Using Hertling's result, we show in Theorem~\ref{thm:localext} that it
is also possible to extend the neighborhood of a smooth point on the
discriminant locus of an $N$-dimensional Frobenius manifold $M$ to a
convergent CohFT, provided that $m = 3$.
Roughly, the CohFT is constructed by deforming the CohFT corresponding
to the $A_2 \times (A_1)^{N-2}$-singularity using the action of an
$R$-matrix constructed from $M$.
We should note that special cases of Theorem~\ref{thm:localext} on the
level of intersection numbers have previously appeared in the
literature, such as \cite{He01} (extension to genus one), \cite{Mi15}
(extension for CohFTs from singularity theory), \cite{CoIr15}
(extension for CohFTs from Gromov--Witten theory).\footnote{It was
  also communicated to the author that Dubrovin--Liu--Zhang have a
  proof for the general conformal case on the level of intersection
  numbers.}
Compared to the existing literature, Theorem~\ref{thm:localext} works
on the cycle level and does not depend on the existence of an Euler
vector field.

\subsection{Tautological relations}

The initial motivation to this work is the study of relations in the
tautological ring of the moduli space of curves.

The tautological rings $RH^*(\Mbar_{g, n})$ are certain subrings of
the cohomology rings $H^*(\Mbar_{g, n})$ of the moduli space
$\Mbar_{g, n}$ of stable curves.
Starting from the 1980s with Mumford's seminal article \cite{Mu83},
they have been studied extensively.
However, their structure is still not completely understood: While
there is an explicit set of generators parameterized by decorated
graphs, the set of relations between the generators is not known.
On the other hand, Pixton's set \cite{Pi12P} of generalized
Faber--Zagier relations gives a well-tested conjectural description
for this set of relations.
Another conjectural description had been given by Faber's Gorenstein
conjecture but it is now known to be false in general \cite{PeTo14}.

In \cite{PPZ15} the relations of Pixton have been shown to arise in
the computation of Witten's 3-spin class, or equivalently the CohFT
corresponding to the $A_2$-singularity, via the Givental--Teleman
reconstruction of semisimple CohFTs.
The formula that Pandharipande--Pixton--Zvonkine obtain for Witten's
3-spin class has the form of a limit $\phi \to 0$ of a Laurent series
in a variable $\phi$ whose coefficients are tautological classes.
The existence of the limit implies cancellation between tautological
classes such that no poles in $\phi$ are left in the end.
These relations between tautological classes, after adding relations
directly following from them, give exactly the relations of Pixton.

As noted in \cite{PPZ15}, the limit $\phi \to 0$ can be viewed as
approaching a non-semisimple point on the Frobenius manifold
corresponding to the $A_2$-singularity.
In particular, the same procedure can be applied to get relations from
other generically semisimple, convergent CohFTs but it is not clear
how the relations from different CohFTs relate to each other.
In \cite{Ja14P} first comparison results have been proven: The
relations from the equivariant Gromov--Witten theory of $\PP^1$ are
equivalent to the relations from the $A_2$-theory and in general the
relations from equivariant $\PP^{N - 1}$ imply the $A_N$-relations.

While the tautological relations obtained from different CohFTs look
very different and are usually more complicated than Pixton's
relations, as an application of Theorem~\ref{thm:localext}, we show in
Theorem~\ref{thm:main} that they give the same set of relations.
This says in some sense that the relations of Pixton are the universal
relations necessary in order for the Givental--Teleman classification
to admit non-semisimple limits.
Theorem~\ref{thm:main} can also be used to relate more geometric
relations to Pixton's relations (see for example \cite{ClJa16P}).

In this paper we work over $\CC$ and with the tautological ring in
cohomology.
It is actually more natural to define the tautological ring in Chow
and the results of this paper apply also to Chow provided that the
Givental--Teleman reconstruction is proven in Chow for the relevant
(Chow-valued) CohFTs.

\subsection*{Plan of the paper}

In Section~\ref{sec:frob} we first recall basic properties of
Frobenius manifolds, and then state Hertling's classification in
Theorem~\ref{thm:struct}.
In Section~\ref{sec:cohft}, we start by recalling the definition of
cohomological field theories and the statement of the
Givental--Teleman classification.
After that, in Section~\ref{sec:cohft:taut}, we discuss the
tautological ring and the relations resulting from the classification.
In Section~\ref{sec:cohft:localext}, we prove
Theorem~\ref{thm:localext} about the local extension of semisimple
Frobenius manifolds.
We discuss in Section \ref{sec:cohft:eq} how its proof implies
Theorem~\ref{thm:main} on the comparison of tautological relations.
In Section~\ref{sec:cohft:glocalext}, we shortly consider the problem
of finding a global extension theorem similar to
Theorem~\ref{thm:localext}.
Afterwards, in Section~\ref{sec:cohft:ex} we study two examples, which
illustrate obstructions to directly generalizing our results.
In the final Section~\ref{sec:cohft:gw} we show that certain other
relations obtained from the equivariant Gromov--Witten theory of toric
targets can also be expressed in terms of Pixton's relations.

\subsection*{Acknowledgments}

The author is very grateful for various discussions with A.~Buryak,
E.~Clader, Y.P.~Lee, D.~Petersen, R.~Pandharipande, C.~Schiessl,
S.~Shadrin and D.~Zvonkine.

The author has learned about the method of obtaining relations by
studying a CohFT near the discriminant from D.~Zvonkine at the
conference \emph{Cohomology of the moduli space of curves} organized
by the \emph{Forschungs-institut für Mathematik at ETH Zürich} in
2013.
D.~Zvonkine there also expressed an idea why we should obtain the same
relations from different CohFTs.
Together with S.~Shadrin he studied them in comparison to relations
obtained from degree considerations (see Remark~\ref{rmk:degree}).

The results on the extension of a Frobenius manifold to a CohFT is
motivated from discussions with M.~Kazarian, T.~Milanov and
D.~Zvonkine at the workshop \emph{Geometric Invariants and Spectral
  Curves} at the \emph{Lorentz Center} in Leiden.

Special thanks are due to an anonymous referee for drawing the
author's attention to Hertling's work on Frobenius manifolds.

This research was carried out while being a PhD. student of
R.~Pandharipande at ETH Zürich and being supported by the Swiss
National Science Foundation grant SNF 200021\_143274.

\section{Frobenius manifolds}
\label{sec:frob}

\subsection{Definition and basic properties}

Frobenius manifolds have been introduced by Dubrovin~\cite{Du94}.
They naturally arise when studying genus zero Gromov--Witten theory.
Let us begin by recalling their basic properties in the following
slightly redundant definition.
\begin{definition}
  An $N$-dimensional (complex, even) Frobenius manifold is a 4-tuple
  $(M, \eta, A, \mathbf 1)$, consisting of
  \begin{itemize}
  \item $M$, a complex, connected manifold of dimension $N$,
  \item a nonsingular metric
    $\eta \in \Gamma(\operatorname{Sym}^2(T^*M))$,
  \item a tensor $A \in \Gamma(\operatorname{Sym}^3(T^*M))$,
  \item a vector field $\mathbf 1 \in \Gamma(TM)$,
  \end{itemize}
  satisfying the following properties:
  \begin{itemize}
  \item A commutative, associative product $\star$ on $TM$, with unit
    $\mathbf 1$, is defined by setting for local vector fields $X$ and
    $Y$ that
    \begin{equation*}
      \eta(X \star Y, Z) = A(X, Y, Z)
    \end{equation*}
    for any local vector field $Z$.
  \item The metric $\eta$ is flat and $\mathbf 1$ is an $\eta$-flat
    vector field.
  \item Locally around each point there exist \emph{flat} coordinates
    $t_{\alpha}$ such that the metric and the unit vector field are
    constant when written in the basis of the corresponding local
    vector fields $\frac \partial{\partial t_\alpha}$.
  \item Locally on $M$ there exists a holomorphic function $\Phi$
    called \emph{potential} such that
    \begin{equation*}
      A\left(\frac \partial{\partial t_\alpha}, \frac \partial{\partial t_\beta}, \frac \partial{\partial t_\gamma}\right) = \frac{\partial^3 \Phi}{\partial t_\alpha\partial t_\beta\partial t_\gamma}.
    \end{equation*}
  \end{itemize}
\end{definition}
\begin{definition}
  \label{def:euler}
  A Frobenius manifold is called \emph{conformal} if it admits an
  \emph{Euler vector field}, i.e.\ a vector field $E$ of the form
  \begin{equation*}
    E = \sum_\mu (\alpha_\mu t_\mu + \beta_\mu) \frac{\partial}{\partial t_\mu},
  \end{equation*}
  such that the quantum product, the unit and the metric are
  eigenfunctions of the Lie derivative $L_E$ with eigenvalues $1$,
  $-1$ and $2 - \delta$ respectively.
  Here $\delta$ is a rational number called \emph{conformal
    dimension}.
\end{definition}

\subsection{Discriminant and semisimplicity}

We say that a Frobenius manifold $M$ is (generically)
\emph{semisimple} if the algebra $(T_p M, \star_p)$ is semisimple for
a generic point $p$ on $M$.
The set of points $p \in M$ such that $(T_p M, \star_p)$ is not
semisimple is called the \emph{discriminant locus} (or the caustic).

If $M$ is semisimple, near any semisimple point, we can choose a basis
$\frac\partial{\partial u_i}$ of orthogonal idempotents and we use the
notation $\Delta_i^{-1}$ for their norms.
Then $\Delta_i^{\frac 12} \frac\partial{\partial u_i}$ define
\emph{normalized idempotents}.
As the notation suggests, the vector fields
$\frac\partial{\partial u_i}$ commute and we can integrate them
locally near semisimple points to give the \emph{canonical
  coordinates} $u_i$.

For this article, the structure of a Frobenius manifold near a smooth
point of the discriminant plays an important role.
We discuss a two-dimensional motivational example before recalling an
important result in this direction due to C.~Hertling.
\begin{example}
  \label{ex:3spin}
  The Givental--Saito theory of the $A_2$-singularity, which appears
  in the study of Witten's 3-spin class, concerns a two-dimensional
  Frobenius manifold.
  As a manifold, it is isomorphic to $\CC^2$ with coordinates
  $t_0, t_1$ and its points correspond to versal transformations
  \begin{equation*}
    \frac{x^3}3 - t_1 x + t_0
  \end{equation*}
  of the $A_2$-singularity $\frac{x^3}3$.
  In the basis $\frac \partial{\partial t_0}$,
  $\frac \partial{\partial t_1}$ the metric $\eta$ is given by the
  matrix
  \begin{equation*}
    \begin{pmatrix}
      0 & 1 \\
      1 & 0
    \end{pmatrix}
  \end{equation*}
  and the potential is
  \begin{equation*}
    \Phi(t_0, t_1) = \frac 12 t_0^2 t_1 + \frac 1{24} t_1^4.
  \end{equation*}
  Therefore $\frac \partial{\partial t_0}$ is the unit, and the only
  interesting quantum product is
  \begin{equation*}
    \frac \partial{\partial t_1} \star \frac \partial{\partial t_1} = t_1 \frac \partial{\partial t_0}.
  \end{equation*}
  Hence on a two-fold cover of $\CC^2$ ramified along the discriminant
  locus $\{t_1 = 0\}$ we can define the meromorphic idempotents
  \begin{equation*}
    \epsilon_\pm = \pm \frac 1{2\sqrt{t_1}} \frac \partial{\partial t_1} + \frac 12 \frac \partial{\partial t_0}.
  \end{equation*}
  A choice of corresponding canonical coordinates is given by
  \begin{equation*}
    u_\pm = t_0 \pm \frac 23 t_1^{\frac 32}.
  \end{equation*}
  Notice that we can recover the flat vector fields by setting
  \begin{equation*}
    \frac \partial{\partial t_0} = \epsilon_+ + \epsilon_-
  \end{equation*}
  and
  \begin{equation*}
    \frac \partial{\partial t_1} = \left(\frac 34(u_+ - u_-)\right)^{\frac 13} (\epsilon_+ - \epsilon_-).
  \end{equation*}
\end{example}

\begin{theorem} (\cite{He99P, He01})
  \label{thm:struct}
  Let $M$ be a (generically) semisimple Frobenius manifold.
  In a neighborhood $U$ of a smooth point of the discriminant of $M$,
  there exists an integer $m \ge 3$ such that
  \begin{itemize}
  \item all but two idempotents $\frac\partial{\partial u_1}$,
    $\frac\partial{\partial u_2}$ extend holomorphically to the
    discriminant in $U$,
  \item for a suitable choice of integration constants, there is a
    holomorphic root $(u_1 - u_2)^{2/m}$, and its vanishing
    locus describes the discriminant in $U$,
  \item The vector fields
    \begin{equation*}
      (u_1 - u_2)^{(m-2)/m}\left(\frac\partial{\partial u_1} - \frac\partial{\partial u_2}\right), \frac\partial{\partial u_1} + \frac\partial{\partial u_2}, \frac\partial{\partial u_{\ge 3}}
    \end{equation*}
    extend holomorphically to the discriminant and span the tangent
    space at every point of $U$.
    The first of these vector fields spans the space of nilpotent
    tangent vectors at $p$.
  \end{itemize}
  Furthermore, if in addition for any flat vector field $X$, for a
  genus one potential $G$, the function $\mathrm dG(X)$ extends to the
  discriminant, the integer $m$ has to be equal to $3$.
\end{theorem}
Here, we recall the formula (see \cite{Gi98}) for a genus one
potential $G$
\begin{equation}
  \label{eq:genus1}
  \mathrm dG = \frac 1{48} \sum_i \mathrm d\log(\Delta_i) + \frac 12 \sum_i r_{ii} \mathrm du_i,  
\end{equation}
in which the functions $r_{ii}$\footnote{They correspond to the
  diagonal entries of the linear part of the $R$-matrix.}, determined
up to an integration constant, satisfy
\begin{equation*}
  \mathrm dr_{ii} = \frac 14 \sum_j \frac{\partial\log(\Delta_j)}{\partial u_i} \frac{\partial\log(\Delta_i)}{\partial u_j} (\mathrm du_j - \mathrm du_i).
\end{equation*}
\begin{proof}[Proof of Theorem~\ref{thm:struct}]
  Semisimple Frobenius manifolds are examples of massive
  $F$-manifolds.
  By Theorem~\cite[Theorem~12.2]{He99P} the germs of a massive
  $F$-manifolds near a smooth point of the discriminant locus is of
  type $I_2(m) \times A_1^{n - 2}$.
  This implies the first part of the Theorem.

  Theorem~\cite[Theorem~6.3]{He01} considers the $G$-function of
  semisimple Frobenius manifolds of type $I_2(m)$, and shows that $m$
  must be equal to 3 when the $G$-function extends to the discriminant
  locus.\footnote{The proof of \cite[Theorem~6.3]{He01} does not
    depend on the existence of an Euler vector field.}
\end{proof}

\section{Cohomological Field Theories}
\label{sec:cohft}

\subsection{Definitions}
\label{sec:cohft:def}

Cohomological field theories were first introduced by Kontsevich and
Manin in \cite{KoMa97} (see also \cite{PPZ15}) to formalize the
structure of classes from Gromov--Witten theory.
Let $V$ be an $N$-dimensional $\CC$-vector space and $\eta$ a
nonsingular bilinear form on $V$.

Before stating their definition, we recall the moduli space
$\Mbar_{g, n}$ of stable, connected, at most nodal algebraic curves of
arithmetic genus $g$ with $n$ markings, as well as that there are
tautological maps induced by forgetting a marking and gluing along
markings:
\begin{align*}
  \Mbar_{g, n + 1} &\to \Mbar_{g, n}, \\
  \Mbar_{g_1, n_1 + 1} \times \Mbar_{g_2, n_2 + 1} &\to \Mbar_{g_1 + g_2, n_1 + n_2}, \\
  \Mbar_{g, n + 2} &\to \Mbar_{g + 1, n}
\end{align*}
\begin{definition}
  A cohomological field theory (CohFT) $\Omega$ on $(V, \eta)$ is a
  system
  \begin{equation*}
    \Omega_{g, n} \in H^*(\Mbar_{g, n}) \otimes (V^*)^{\otimes n}
  \end{equation*}
  of multilinear forms with values in the cohomology ring of $\Mbar_{g, n}$
  satisfying the following properties:
  \begin{description}
  \item[$S_n$-Equivariance] $\Omega_{g, n}$ is $S_n$-equivariant with
    respect to the $S_n$-action permuting the markings of
    $\Mbar_{g, n}$ and the factors of $(V^*)^{\otimes n}$.
  \item[Gluing] The pull-back of $\Omega_{g, n}$ via the gluing map
    \begin{equation*}
      \Mbar_{g_1, n_1 + 1} \times \Mbar_{g_2, n_2 + 1} \to \Mbar_{g, n}
    \end{equation*}
    is given by the direct product of $\Omega_{g_1, n_2 + 1}$ and
    $\Omega_{g_2, n_2 + 1}$ with the bivector $\eta^{-1}$ inserted at
    the two points glued together.
    Similarly, for the gluing map
    $\Mbar_{g - 1, n + 2} \to \Mbar_{g, n}$, the pull-back of
    $\Omega_{g, n}$ is given by $\Omega_{g - 1, n + 2}$ with
    $\eta^{-1}$ inserted at the two points glued together.
  \item[Unit] There is a special element $\mathbf 1 \in V$ called the
    \emph{unit} such that
    \begin{equation*}
      \Omega_{g, n + 1}(v_1, \dotsc, v_n, \mathbf 1)
    \end{equation*}
    is the pull-back of $\Omega_{g, n}(v_1, \dotsc, v_n)$ under the
    forgetful map and
    \begin{equation*}
      \Omega_{0, 3}(v, w, \mathbf 1) = \eta(v, w).
    \end{equation*}
  \end{description}
\end{definition}
\begin{definition}
  A CohFT which takes values in $H^0$ (so that each $\Omega_{g, n}$ is
  a multiple of the fundamental class) is called a \emph{topological
    field theory} (TQFT).
\end{definition}

The definition of CohFTs can be also generalized to families of CohFTs
over a ground ring.
We will use the following non-standard definition.
\begin{definition}
  Let $e_1, \dotsc, e_N$ be a basis of $V$.
  A \emph{convergent CohFT} $\Omega$ on $V$ is a CohFT defined over
  the ring of holomorphic functions of an open neighborhood $U$ of
  $0 \in V$ such that for all $g \ge 0$, all
  $\alpha_1, \dotsc, \alpha_n \in V$ and all
  $\mathbf t = t_1 e_1 + \dotsb + t_N e_N \in U$ we have
  \begin{equation}
    \label{eq:shift}
    \Omega_{g, n}|_{\mathbf t}(\alpha_1, \dotsc, \alpha_n) = \sum_{k = 0}^\infty \frac 1{k!} \pi_* \Omega_{g, n + k}|_0(\alpha_1, \dotsc, \alpha_n, \mathbf t, \dotsc, \mathbf t),
  \end{equation}
  where $\pi\colon \Mbar_{g, n + k} \to \Mbar_{g, n}$ forgets the last
  $k$ markings.
\end{definition}
\begin{remark}
  If $\Omega$ is defined via the Gromov--Witten theory of $X$, then
  convergence of $\Omega$ implies convergence of the total ancestor
  potential of $X$ in the sense of \cite{CoIr15}.
  This follows immediately by integration of $\Omega_{g, n}$ over
  $\Mbar_{g, n}$ against monomials in $\psi$-classes.
  A priori, the integrals of $\Omega$ against $\psi$-classes contain
  less information than $\Omega$ itself, so it is not clear whether
  the convergence of the CohFT is equivalent to the convergence of its
  ancestor potential.
  Still, it should be straightforward to adapt many proofs of
  convergence of the ancestor potential to show convergence of the
  corresponding CohFT level (see also Section~\ref{sec:cohft:gw}).
\end{remark}
We can define from any usual CohFT a convergent CohFT by using
\eqref{eq:shift}, under the assumption that the sum converges in a
neighborhood of 0.
\begin{definition}
  The \emph{underlying Frobenius manifold} of a convergent CohFT
  $\Omega$ is, as a manifold, the neighborhood $U$ of $0 \in V$.
  At every point $p$ of $U$, the tangent space is identified with $V$
  via the map sending $\frac\partial{\partial t_\mu}$ at $p$ to
  $e_\mu$.
  With this identification, $\eta$ defines the metric, $\mathbf 1$
  defines the unit vector field and $\Omega_{0, 3}$ defines the
  symmetric tensor $A$.
\end{definition}
\begin{remark}
  Restricting to the origin, we see that every CohFT determines a
  Frobenius algebra.
  This operation restricts to a bijection between TQFTs and Frobenius
  algebras of dimension $N$.
\end{remark}

Using the underlying Frobenius manifold, for any convergent CohFT we
can define the quantum product on $V$ (depending on a point in $U$),
semisimplicity and the discriminant locus.

\begin{example}
  \label{ex:extend}
  Given an $N$-dimensional (convergent) CohFT $\Omega$ and some
  $c \in \CC^*$, we can define an $(N + 1)$-dimensional (convergent)
  CohFT $\Omega'$: If $V$ is the underlying vector space of $\Omega$,
  then $V \oplus \langle v\rangle$ will be the underlying vector space
  of $\Omega'$.
  The nonsingular bilinear form $\eta'$ on $V \oplus \CC$ is defined
  via $\eta'(v, v) = c$, $\eta'(\alpha, v) = 0$ and
  $\eta'(\alpha, \beta) = \eta(\alpha, \beta)$, where
  $\alpha, \beta \in V$ and $\eta$ is the nonsingular bilinear form of
  $V$.
  The CohFT $\Omega'$ is then defined by multilinearity from setting
  \begin{equation*}
    \Omega'_{g, n}(\alpha_1, \dotsc, \alpha_n) = \Omega_{g, n}(\alpha_1, \dotsc, \alpha_n),
  \end{equation*}
  if all $\alpha_i$ lie in $V$, imposing the condition that $\Omega'$
  vanishes if one argument is a multiple of $v$ and another argument
  lies in $V$, and setting
  \begin{equation*}
    \Omega'_{g, n}(v, \dotsc, v) = c^{1 - g}.
  \end{equation*}
  Finding the right definition in the remaining case $n = 0$ and
  checking the axioms of a CohFT is left as an exercise to the reader.

  Notice that $v$ will be an idempotent element for the quantum
  product and that this operation therefore preserves semisimplicity.
\end{example}

\begin{definition}
  A convergent CohFT $\Omega_{g, n}$ is called \emph{homogeneous} if its
  underlying Frobenius manifold is conformal and the CohFT is an
  eigenvector of $L_E$ of eigenvalue $(g - 1)\delta + n$.
\end{definition}
As the name suggests, a convergent CohFT is homogeneous if it carries a
grading such that all natural structures are homogeneous with respect
to the grading.

\subsection{Reconstruction}
\label{sec:cohft:reconstr}

The (upper half of the) \emph{symplectic loop group} corresponding to
a vector space $V$ with nonsingular bilinear form $\eta$ is the group
of endomorphism valued power series $V[\![z]\!]$ such that the
\emph{symplectic condition} $R(z) R^t(-z) = 1$ holds.
Here $R^t$ is the adjoint of $R$ with respect to $\eta$.
There is an action of this group on the space of all CohFTs based on a
fixed semisimple Frobenius algebra structure of $V$.
The action is often named after Givental because he has introduced it
on the level of arbitrary genus Gromov--Witten potentials.

Given a CohFT $\Omega_{g, n}$ and such an endomorphism $R$, the new
CohFT $R\Omega_{g, n}$ takes the form of a sum over dual graphs
$\Gamma$
\begin{equation}
  \label{eq:reconstra}
  R\Omega_{g, n}(v_1, \dotsc, v_n)
  = \sum_\Gamma \frac 1{\Aut(\Gamma)} \xi_*\left(\prod_v \sum_{k = 0}^\infty \frac 1{k!} \pi_* \Omega_{g_v, n_v + k}(\dots)\right),
\end{equation}
where $\xi\colon \prod_v \Mbar_{g_v, n_v} \to \Mbar_{g, n}$ is the
gluing map of curves of topological type $\Gamma$ from their
irreducible components,
$\pi\colon \Mbar_{g_v, n_v + k} \to \Mbar_{g_v, n_v}$ forgets the last
$k$ markings and we still need to specify what is put into the
arguments of $\prod_v \Omega_{g_v, n_v + k}$.
Instead of only allowing vectors in $V$ to be put into $\Omega_{g, n}$
we will also allow for elements of $V[\![\psi_1, \dotsc, \psi_n]\!]$
where $\psi_i$ acts on the cohomology of the moduli space by
multiplication with the $i$th cotangent line class.
\begin{itemize}
\item Into each argument corresponding to a marking of the curve, put
  $R^{-1}(\psi)$ applied to the corresponding vector.
\item Into each pair of arguments corresponding to an edge put the
  bivector
  \begin{equation*}
    \frac{R^{-1}(\psi_1) \otimes R^{-1}(\psi_2) - \Id}{-\psi_1 - \psi_2} \eta^{-1} \in V^{\otimes 2}[\![\psi_1, \psi_2]\!],
  \end{equation*}
  where one has to substitute the $\psi$-classes at each side of the
  normalization of the node for $\psi_1$ and $\psi_2$.
  By the symplectic condition this is well-defined.
\item At each of the additional arguments for each vertex put
  \begin{equation*}
    T(\psi) := \psi(\Id - R^{-1}(\psi)) \mathbf 1,
  \end{equation*}
  where $\psi$ is the cotangent line class corresponding to that
  vertex.
  Since $T(z) = O(z^2)$ the above $k$-sum is finite.
\end{itemize}

The following reconstruction result (on the level of potentials) has
been first proposed by Givental \cite{Gi01b}.
\begin{theorem}[\cite{Te12}]
  \label{thm:reconstruction}
  The $R$-matrix action is free and transitive on the space of
  semisimple CohFTs based on a given Frobenius algebra.

  Furthermore, given a convergent semisimple CohFT $\Omega$, locally
  around a semisimple point, the element $R$ of the symplectic loop
  group, taking the TQFT corresponding to the Frobenius algebra to
  $\Omega$, satisfies the following differential equation of one-forms
  when written in a basis of normalized idempotents
  \begin{equation}
    \label{eq:RDE}
    [R(z), d\mathbf u] + z\Psi^{-1} d(\Psi R(z)) = 0.
  \end{equation}
  Here $\mathbf u$ is the diagonal matrix filled with the canonical
  coordinates $u_i$ corresponding to the idempotents and $\Psi$ is the
  basis change from the basis of normalized idempotents to a flat
  basis.
\end{theorem}
\begin{remark}
  The differential equation \eqref{eq:RDE} makes sense for any
  Frobenius manifold.
  In general it defines $R$ only up to right multiplication by a
  diagonal matrix whose entries are of the form
  $\exp(a_1 z + a_3z^3 + \dotsb)$, where the $a_i$ are constant on the
  Frobenius manifold \cite{Gi01a}.
  If the further condition of homogeneity with respect to an Euler
  vector field is imposed on $R$, there is a unique solution.
\end{remark}
\begin{remark}
  Teleman's proof relies heavily on topological results (Mumford's
  conjecture/Madsen--Weiss theorem) and it is therefore not known if
  the same classification result also holds in general when we work in
  Chow instead of cohomology.
  It is still known that the Chow-valued statement holds in some cases
  such as for the equivariant Gromov--Witten theory of a toric
  variety.
\end{remark}
\begin{remark}
  Formula~\eqref{eq:genus1} for a genus one potential
  \begin{equation*}
    \mathrm dG(X) = \int_{\Mbar_{1, 1}} \Omega_{1, 1}(X)
  \end{equation*}
  is a special case of the reconstruction.
\end{remark}

For later use, let us use the normalized idempotent basis to make the
reconstruction formula a bit more concrete:
We can decompose any endomorphism $F$ of $V$ into a collection of
linear forms
\begin{equation*}
  F = \sum_i F^i \tilde\epsilon_i,
\end{equation*}
where $\tilde\epsilon_i$ is the $i$th normalized idempotent element,
and we will use the formula
\begin{equation*}
  \omega_{g, n} (\tilde\epsilon_{a_1}, \dotsc, \tilde\epsilon_{a_n}) =
  \begin{cases}
    \sum_i \Delta_i^{g - 1}, & \text{if } n = 0, \\
    \Delta_{a_1}^{\frac{2g - 2 + n}2}, & \text{if } a_1 = \dotsb = a_n, \\
    0, & \text{else,}
  \end{cases}
\end{equation*}
where the $\Delta_i$ are the inverses of the norms of the idempotents.
Then we can rewrite \eqref{eq:reconstra} to
\begin{equation}
  \label{eq:reconstrb}
  R\Omega_{g, n}(v_1, \dotsc, v_n)
  = \sum_{\Gamma, c} \frac 1{\Aut(\Gamma, c)} \xi_*\left(\prod_v C_{v, c(v)}(\dotsc)\right),
\end{equation}
where $c$ is a coloring of the vertices of $\Gamma$ by a color in the
set $\{1, \dotsc, N\}$, and the local contribution $C_{v, i}$ at a
vertex $v$ of genus $g$, with $n$ markings and of color $i$ is an
$n$-form taking power series in $z$ as inputs, and is given by
\begin{multline*}
  C_{v, i}(\alpha_1, \dotsc, \alpha_n) \\
  = \sum_{k = 0}^\infty \frac{\Delta_i^{\frac{2g - 2 + n + k}2}}{k!} \pi_*\Bigg(\prod_{j = 1}^n \alpha_j(\psi_j) \prod_{j = n + 1}^{n + k} \psi_j (\Id^i - (R^{-1}(\psi_j))^i) \mathbf 1_\Omega\Bigg).
\end{multline*}
The still missing arguments in \eqref{eq:reconstrb}, which correspond
to preimages of the marked points and nodes in the normalization, are
the same as the vectors and bivectors also used in
\eqref{eq:reconstra}, except that we need to take their coordinates
as prescribed by the coloring.

\subsection{Relations in the tautological ring}
\label{sec:cohft:taut}

The tautological subrings $R^*(\Mbar_{g, n})$ can be compactly defined
\cite{FaPa05} as the smallest system of subrings
$R^*(\Mbar_{g, n}) \subseteq H^*(\Mbar_{g, n})$ stable under
push-forward under the tautological maps as described in
Section~\ref{sec:cohft:def}.
Each tautological ring has a finite additive generating set
\cite{GrPa03} which is indexed by strata of $\Mbar_{g, n}$
(corresponding to dual graphs) decorated with Morita--Mumford--Miller
$\kappa$-classes and $\psi$-classes.
The formal algebra generated by these generators is called the strata
algebra $\mathcal S_{g, n}$ \cite{Pi12P}.
Pushforwards and pullbacks along the gluing and forgetful morphisms
can be lifted to homomorphisms of the corresponding strata algebras
satisfying the push-pull formula, \dots. 
Relations in the tautological ring are elements of the kernel of the
natural projection $q\colon \mathcal S_{g, n} \to R^*(\Mbar_{g, n})$.

Consider a semisimple, $N$-dimensional convergent CohFT $\Omega$
defined in a neighborhood $U$ of $0 \in V$.
Let $D \subset U$ be the discriminant locus.
By the reconstruction formula described in
Section~\ref{sec:cohft:reconstr}, for each point outside $D$ in $U$ we
can find an $R$-matrix such that $\Omega$ is given by applying the
action of $R$ to the underlying TQFT.

We obtain relations in the tautological ring by studying the behavior
along $D$.
On the one hand, the reconstruction formula defines a function
$\widetilde\Omega_{g, n} \in \mathcal S_{g, n} \otimes (V^*)^{\otimes
  n} \otimes \cO_{U \setminus D}$ that might have singularities along
the discriminant locus.\footnote{The proof of Theorem~\ref{thm:main}
  will imply that for a CohFT the only possible singularities of the
  $R$-matrix along $D$ are poles.}
On the other hand, by Theorem~\ref{thm:reconstruction}, we recover the
original CohFT if we project from the strata algebra to the
tautological ring, so that
$q(\widetilde\Omega_{g, n}) \in R^*(\Mbar_{g, n}) \otimes
(V^*)^{\otimes n} \otimes \cO_U$.
Therefore, if
$p\colon \cO_{U \setminus D} \to (\cO_{U \setminus D}/\cO_U)$ denotes
the projection, then
$p(\widetilde\Omega_{g, n}) \in \ker(q) \otimes (V^*)^{\otimes n}
\otimes (\cO_{U \setminus D}/\cO_U)$ is a
$(V^*)^{\otimes n} \otimes (\cO_{U \setminus D}/\cO_U)$-valued
tautological relation, from which we can extract many ordinary
tautological relations by applying linear forms
$\varphi\colon (V^*)^{\otimes n} \otimes (\cO_{U \setminus D}/\cO_U)
\to \CC$.
\begin{definition}
  \label{def:rels}
  The vector space of tautological relations associated to the
  convergent CohFT $\Omega$ is defined as the smallest system of
  ideals of $\mathcal S_{g, n}$ which is stable under push-forwards
  via the gluing and forgetful morphisms and contains the relations
  $p(\widetilde\Omega_{g, n})$ from cancellations of singularities in
  the reconstruction of $\Omega$, that we have just discussed.
\end{definition}
\begin{example}
  For the 2-dimensional, (convergent) CohFT corresponding to Witten's
  3-spin class, in \cite{PPZ15} it is proven that the ideal of
  relations coincides with the relations of Pixton \cite{Pi12P}, which
  are conjectured to be all relations between tautological classes.
\end{example}
\begin{example}
  In \cite{Ja14P} it is shown that the ideal of relations of the
  Gromov--Witten theory of equivariant projective space $\PP^{N - 1}$
  contains the relations for Witten's $(N + 1)$-spin class.
\end{example}
\begin{example}
  In \cite{Ja14P} it is also shown that the set of relations for
  equivariant $\PP^1$ and Witten's 3-spin class coincide.
\end{example}
\begin{remark}
  For nonequivariant $\PP^1$, the theory does not apply since the
  Frobenius manifold is semisimple at all points.
  There is however a different way of how to extract relations in this
  case, which we will study in Section~\ref{sec:cohft:gw}.
\end{remark}

The following is our main result.
\begin{theorem}
  \label{thm:main}
  For any two semisimple convergent CohFTs which are not semisimple at
  all points of the underlying Frobenius manifold, the sets of
  associated tautological relations coincide.
\end{theorem}
\begin{remark}
  In the proof of Theorem~\ref{thm:main} we will first, locally near a
  smooth point on the discriminant, identify canonical coordinates and
  normalized idempotents.
  An important part of the proof is to show that under this
  identification the quotient of corresponding $R$-matrices is
  holomorphic along the discriminant.

  In \cite{Ja14P} for the comparison of equivariant $\PP^1$ and the
  $A_2$-singularity a different, more explicit identification of
  coordinates and vector fields is chosen.
  Therefore, while with this identification the quotient of the
  $R$-matrices is not holomorphic along the discriminant, there exists
  a holomorphic function $\varphi$ such that
  $R_{\PP^1}(z) = R(z) R_{A_2}(\varphi z)$.
  This result depends on the fact that the $A_2$-theory is
  homogeneous.
\end{remark}
\begin{remark}
  \label{rmk:degree}
  The more classical way of \cite{PPZ15} and \cite{PPZ16P} to obtain
  tautological relations from a convergent CohFT works by considering
  cohomological degrees: Assume that $\Omega$ is in addition
  homogenous for an Euler vector field $E$, that all $\beta_i$ vanish
  and that all $\alpha_i$ as in Definition~\ref{def:euler} are
  positive.
  Then the homogeneity implies that the cohomological degree of
  $\Omega_{g, n}(\frac\partial{\partial t_{a_1}}, \dotsc,
  \frac\partial{\partial t_{a_n}})$ is bounded by
  \begin{equation}
    \label{eq:relbound}
    (g - 1)\delta + n - \sum_j \alpha_{a_j}.
  \end{equation}
  However the reconstructed theory might also contain terms of higher
  cohomological degree.
  These thus have to vanish, giving tautological relations.
  To make them more concrete, the coefficients of these relations,
  which are functions on the underlying Frobenius manifold, can be
  evaluated at any point outside of the discriminat.

  These relations from degree considerations follow from the relations
  of Definition~\ref{def:rels}.
  To see this, notice that by positivity of the $\alpha_i$, the action
  of $E$ on holomorphic functions on the underlying Frobenius manifold
  has only positive eigenvalues.
  Therefore and because of homogeneity (and generic semisimplicity),
  coefficients of the part of the reconstruction formula violating the
  bound \eqref{eq:relbound} in particular cannot be holomorphic along
  the discriminant.
\end{remark}

\subsection{Local extension}
\label{sec:cohft:localext}

The proof of the following theorem will occupy this section.
The content of the proof is also used for proving
Theorem~\ref{thm:main}.
\begin{theorem}
  \label{thm:localext}
  Let $M$ be an $N$-dimensional semisimple Frobenius manifold and let
  $p$ be a smooth point of the discriminant of $M$ such that $m = 3$
  in Theorem~\ref{thm:struct}.
  Then, after possibly shrinking $M$ to a smaller neighborhood of $p$,
  there exists a convergent CohFT with underlying Frobenius manifold
  $M$.
\end{theorem}

We first study the consequences of Theorem~\ref{thm:struct} in more
detail.
After possibly shrinking $M$, it gives us a basis of holomorphic
vector fields
$\{\frac\partial{\partial t_0}, \frac\partial{\partial t},
\frac\partial{\partial u_{\ge 3}}\}$, where
\begin{equation}
  \label{eq:basiscohft}
  \frac\partial{\partial t_0} = \frac\partial{\partial u_1} + \frac\partial{\partial u_2}, \qquad \frac\partial{\partial t} = \left(\frac 34(u_1 - u_2)\right)^{\frac 13} \left(\frac\partial{\partial u_1} - \frac\partial{\partial u_2}\right).
\end{equation}
It is easy to see that these vector fields commute and therefore we
can integrate them to coordinates $t_0$, $t$ and $u_{\ge 3}$.
The discriminant locus $D$ is then locally given by the equation
$t = 0$.

Notice that there is a root $\sqrt{t}$ of $t$ such that
\begin{equation*}
  \frac\partial{\partial t} = \sqrt{t} \left(\frac\partial{\partial u_1} - \frac\partial{\partial u_2}\right), \qquad
  \left(\frac\partial{\partial t}\right)^2 = t \left(\frac\partial{\partial u_1} + \frac\partial{\partial u_2}\right).
\end{equation*}
Define holomorphic functions $\eta_0$ and $\eta_1$ by
\begin{equation*}
  \eta_0 = \eta\left(\frac\partial{\partial t_0}, \frac\partial{\partial t_0}\right), \qquad
  \eta_1 = \eta\left(\frac\partial{\partial t}, \frac\partial{\partial t_0}\right)
\end{equation*}
and notice that
\begin{equation*}
  \eta\left(\frac\partial{\partial t}, \frac\partial{\partial t}\right) = \eta\left(\frac\partial{\partial t} \star \frac\partial{\partial t}, \sum_{i = 1}^N \frac\partial{\partial u_i}\right) = t \eta_0.
\end{equation*}
Since $\eta$ is nonsingular, $\eta_1$ cannot vanish on the
discriminant.
The inverses $\Delta_1$, $\Delta_2$ of the norms of the first two
idempotents are given by
\begin{equation*}
  \Delta_1 = \frac{2\sqrt t}{\eta_1 + \sqrt t \eta_0}, \qquad \Delta_2 = \frac{-2\sqrt t}{\eta_1 - \sqrt t \eta_0}.
\end{equation*}

We next choose roots $\sqrt{2 \sqrt t}$, $\sqrt{-2 \sqrt t}$ and
$\sqrt{\eta_1}$.
These induce roots of $\Delta_1$, $\Delta_2$, which we will use to
define the normalized idempotents.
Let $\Psi_0$ be the block diagonal matrix with upper left block being
\begin{equation}
  \label{eq:Psi0}
  \begin{pmatrix}
    \frac{\sqrt t}{\sqrt{2 \sqrt t}} & \frac{-\sqrt t}{\sqrt{-2 \sqrt t}} \\
    \frac 1{\sqrt{2 \sqrt t}} & \frac 1{\sqrt{-2 \sqrt t}}
  \end{pmatrix}
\end{equation}
and the identity matrix as the lower right block.
For $\Psi_0^{-1}$ the upper left block is given by
\begin{equation}
  \label{eq:Psi0inv}
  \begin{pmatrix}
    \frac 1{\sqrt{2 \sqrt t}} & \frac{\sqrt t}{\sqrt{2 \sqrt t}} \\
    \frac 1{\sqrt{-2 \sqrt t}} & \frac{-\sqrt t}{\sqrt{-2 \sqrt t}}
  \end{pmatrix}.
\end{equation}
The matrix $\Psi_0$ is the basis change from normalized idempotents to
the basis
$\{\frac\partial{\partial t_0}, \frac\partial{\partial t},
\frac\partial{\partial u_{\ge 3}}\}$.
In the $A_2$-singularity case $\eta_0 = 0$, $\sqrt{\eta_1} = 1$ and
$\sqrt{\Delta_{\ge 3}} = 1$.

Let $\Psi_1$ denote the basis change from the normalized idempotent
basis to a flat basis and define $\tilde \Psi_1 = \Psi_1 \Psi_0^{-1}$.
\begin{lemma}
  \label{lem:basischange}
  The basis change matrix $\tilde\Psi_1$ is holomorphic along $D$.
\end{lemma}
\begin{proof}
  By Theorem~\ref{thm:struct} it is enough to prove the same statement
  for $\tilde \Psi' := \Psi' \Psi_0^{-1}$ where $\Psi'$ is the basis
  change from the normalized idempotent basis to the basis
  $\{\frac\partial{\partial t_0}, \frac\partial{\partial t},
  \frac\partial{\partial u_{\ge 3}}\}$.
  Since the basis changes leave all but the first two idempotents
  invariant we will only need to consider the upper-left $2 \times 2$
  block of $\tilde\Psi'$.
  We factor this block into
  \begin{equation*}
    \begin{pmatrix}
      \frac{\sqrt t}{2 \sqrt t} & \frac{-\sqrt t}{-2 \sqrt t} \\
      \frac 1{2 \sqrt t} & \frac 1{-2 \sqrt t}
    \end{pmatrix}
    \begin{pmatrix}
      \frac{\sqrt{\Delta_1}}{\sqrt{2t}} & 0 \\
      0 & \frac{\sqrt{\Delta_2}}{\sqrt{-2t}}
    \end{pmatrix}
    \begin{pmatrix}
    1 & \sqrt t \\
    1 & -\sqrt t
  \end{pmatrix},
  \end{equation*}
  a change from
  $\{\frac\partial{\partial t_0}, \frac\partial{\partial t}\}$ to the
  idempotents, a multiplication by a diagonal matrix and the change
  back from the idempotents to
  $\{\frac\partial{\partial t_0}, \frac\partial{\partial t}\}$.
  So we see that the upper left block of $\tilde\Psi'$ has the form
  \begin{equation*}
    \eta_1^{-\frac 12}
    \begin{pmatrix}
      a & t c \\
      c & a
    \end{pmatrix},
  \end{equation*}
  where
  \begin{align*}
    a = \frac 12 \left(\frac 1{\sqrt{1 + \sqrt t \frac{\eta_0}{\eta_1}}} + \frac 1{\sqrt{1 - \sqrt t \frac{\eta_0}{\eta_1}}}\right) =& 1 + \frac 38 \frac{\eta_0^2}{\eta_1^2} t + O(t^2) \\
    c = \frac 1{2\sqrt t} \left(\frac 1{\sqrt{1 + \sqrt t \frac{\eta_0}{\eta_1}}} - \frac 1{\sqrt{1 - \sqrt t \frac{\eta_0}{\eta_1}}}\right) =& -\frac 12 \frac{\eta_0}{\eta_1} - \frac 5{16} \frac{\eta_0^3}{\eta_1^3} t + O(t^2)
  \end{align*}
  which are holomorphic along the discriminant.
\end{proof}

We can repeat this setup for any other $N$-dimensional Frobenius
manifold satisfying the assumptions of Theorem~\ref{thm:localext}, in
particular, as we will assume from now on, for any Frobenius manifold
$M^2$ underlying an $N$-dimensional convergent CohFT $\Omega^2$, such
as an extension of the theory of the $A_2$-singularity to $N$
dimensions using the construction of Example~\ref{ex:extend}
repeatedly.
Our main strategy in order to prove Theorem~\ref{thm:localext} is to
show that a local extension $\Omega^1$ can be obtained from $\Omega^2$
by the action of a holomorphic $R$-matrix and a holomorphic shift.

For this, we first use the coordinates $t, t_0, u_{\ge 3}$ to identify
a small neighborhood of $p$ with a small neighborhood of the origin of
$M^2$.
Let us shrink $M$ accordingly.
Notice that this isomorphism of complex manifolds, if the third roots
in \eqref{eq:basiscohft} have been chosen compatibly, amounts, outside
of the discriminant, to identifying their canonical coordinates.
We accordingly identify normalized idempotents and thereby the tangent
spaces that they span.
This identification preserves the metric but not the quantum product
structure.
In particular, the basis change $\Psi_2$ from the normalized
idempotent basis to a basis of flat vector fields on $M^2$ in general
does not agree with $\Psi_1$.
We set $\tilde \Psi_2 = \Psi_2 \Psi_0^{-1}$, which by
Lemma~\ref{lem:basischange} is also holomorphic along $D$.

The following lemma constructs the holomorphic $R$-matrix we will use
in the construction of $\Omega^1$ from $\Omega^2$.
\begin{lemma}
  \label{lem:Rreg}
  There exists a symplectic solution $R_1$ of the flatness equation
  \eqref{eq:RDE} for $\Psi_1$ such that if $R_2$ denotes the solution
  of \eqref{eq:RDE} for $\Psi_2$ used for reconstructing the CohFT
  $\Omega^2$, the endomorphism $R_1 R_2^{-1}$ is holomorphic (under
  the identifications we have made above).
\end{lemma}
\begin{proof}
  For this proof let $R_1$ and $R_2$ denote the $R$-matrices written
  in the basis of normalized idempotents instead of the underlying
  endomorphism valued power series.
  We set
  \begin{equation*}
    R_i = \Psi_0^{-1} \tilde R_i \Psi_0,
  \end{equation*}
  where $\Psi_0$ is as in \eqref{eq:Psi0}.
  We can write the flatness equations \eqref{eq:RDE} as
  \begin{equation*}
    [\tilde R_i, \Psi_0 \mathrm du \Psi_0^{-1}] + z \tilde\Psi_i^{-1} \mathrm d(\tilde\Psi_i \tilde R_i \Psi_0)\Psi_0^{-1} = 0.
  \end{equation*}
  If $R := \tilde R_1 \tilde R_2^{-1}$, these two differential
  equations combine to
  \begin{equation}
    \label{eq:mixed}
    0 = [R, \Psi_0 \mathrm du \Psi_0^{-1}] + z \tilde\Psi_1^{-1} \mathrm d(\tilde\Psi_1 R \tilde\Psi_2^{-1}) \tilde\Psi_2.
  \end{equation}
  By Lemma~\ref{lem:basischange} it is enough to show that there is a
  solution $R$ of \eqref{eq:mixed} all of whose entries are
  holomorphic along the discriminant and which satisfies the
  symplectic condition.

  We analyze the entries of the ingredients in \eqref{eq:mixed}.
  For this we will consider all matrices to consist of four blocks
  numbered according to
  \begin{equation*}
    \begin{pmatrix}
      1. & 2. \\
      3. & 4.
    \end{pmatrix},
  \end{equation*}
  such that the first block has size $2 \times 2$.
  By Lemma~\ref{lem:basischange} the matrices $\tilde\Psi_i^{-1}$,
  $\tilde\Psi_i$ for $i \in \{1, 2\}$ and therefore also the matrices
  $\tilde\Psi_1^{-1} \mathrm d\tilde\Psi_1$ and
  $(\mathrm d\tilde\Psi_2^{-1}) \tilde\Psi_2$ of one-forms are
  holomorphic along $D$.
  The matrix $(\mathrm d\Psi_0) \Psi_0^{-1}$ has all blocks equal to
  zero except for the first one, which is
  \begin{equation}
    \label{eq:dPsi0}
    \frac{\mathrm dt}{4t}
    \begin{pmatrix}
      1 & 0 \\
      0 & -1
    \end{pmatrix}
  \end{equation}
  and the matrix $\Psi_0 \mathrm d\mathbf u \Psi_0^{-1}$ is block
  diagonal with first block being
  \begin{equation}
    \label{eq:PsiduPsi}
    \begin{pmatrix}
      \mathrm dt_0 & t \mathrm dt \\
      \mathrm dt & \mathrm dt_0
    \end{pmatrix}
  \end{equation}
  and the other block being the diagonal matrix with entries
  $\mathrm du_{\ge 3}$.
  Furthermore, because in general $\Psi_1^{-1} \mathrm d\Psi_1$ is
  antisymmetric (as can be seen by differentiating
  $\Psi_1^t \eta \Psi_1 = 1$) and from \eqref{eq:Psi0},
  \eqref{eq:Psi0inv} and \eqref{eq:dPsi0} we see that the first block
  of $\tilde\Psi_1^{-1} \mathrm d\tilde\Psi_1$ has the general form
  \begin{equation}
    \label{eq:dPsitilde}
    \begin{pmatrix}
      x & 0 \\
      0 & -x
    \end{pmatrix}
  \end{equation}
  and the forth block is still antisymmetric.
  The same holds for $(\mathrm d\tilde\Psi_2^{-1}) \tilde\Psi_2$.

  We construct the coefficients of $R$ inductively.
  Let us set
  \begin{equation*}
    R(z) = \sum_{i = 0}^\infty R^i z^i
  \end{equation*}
  and $R_{jk}^i$ for the entries of $R^i$.
  We assume that we have already constructed $R^j$ for $j \le i$
  satisfying the flatness equation and symplectic condition modulo
  $z^{i + 1}$.

  Because of \eqref{eq:PsiduPsi}, inserting $\frac\partial{\partial
    t_0}$ into the $z^{i + 1}$-part of \eqref{eq:mixed} directly gives
  equations for the off-diagonal blocks of $R^{i + 1}$ in terms of
  holomorphic functions.

  Similarly, inserting $\frac\partial{\partial u_{\ge 3}}$ into the
  $z^{i + 1}$-part of \eqref{eq:mixed} gives us holomorphic formulas
  for the off-diagonal entries of $R^{i + 1}$ in the forth block.
  For the diagonal entries of this block we instead insert
  $\frac\partial{\partial t}$ into the $z^{i + 2}$-part of
  \eqref{eq:mixed} and because of the antisymmetry obtain that the
  first $t$-derivatives of the diagonal entries are holomorphic.
  We can integrate them locally and have an arbitrary choice of
  integration constants (ignoring the symplectic condition for the
  moment).

  It remains the analysis of the first block.
  For this it is useful to compute the commutator
  \begin{equation*}
    \left[
    \begin{pmatrix}
      R_{11}^{i + 1} & R_{12}^{i + 1} \\
      R_{21}^{i + 1} & R_{22}^{i + 1}
    \end{pmatrix}
    ,
    \begin{pmatrix}
      0 & t \\
      1 & 0
    \end{pmatrix}
    \right]
    =
    \begin{pmatrix}
      R_{12}^{i + 1} - tR_{21}^{i + 1} & t(R_{11}^{i + 1} - R_{22}^{i + 1}) \\
      R_{22}^{i + 1} - R_{11}^{i + 1} & tR_{21}^{i + 1} - R_{12}^{i + 1}.
    \end{pmatrix}
  \end{equation*}
  So the insertion of $\frac\partial{\partial t}$ into the
  $z^{i + 1}$-part of \eqref{eq:mixed} gives holomorphic formulas for
  $R_{11}^{i + 1} - R_{22}^{i + 1}$ and
  $R_{12}^{i + 1} - tR_{21}^{i + 1}$.
  Because of \eqref{eq:dPsi0} and \eqref{eq:dPsitilde} the
  $\frac\partial{\partial t}$ insertion into the $z^{i + 2}$-part of
  \eqref{eq:mixed} shows that the first $t$-derivative of
  $R_{11}^{i + 1} + R_{22}^{i + 1}$ is holomorphic and therefore this
  sum is holomorphic and we again have the choice of an integration
  constant.
  Similarly, we find that
  \begin{equation*}
    2t\frac\partial{\partial t} R_{21}^{i + 1} + R_{21}^{i + 1}.
  \end{equation*}
  is holomorphic in $D$ and therefore $R_{21}^{i + 1}$ is holomorphic
  up to a possible constant multiple of $t^{-\frac 12}$.
  Here, we have a unique choice of integration constant giving a
  holomorphic solution.

  In general the symplectic condition does not constrain the
  integration constants of $R^{i + 1}$ when $i + 1$ is odd
  \cite{Gi01a}.
  On the other hand, it completely determines the integration
  constants of $R^{i + 1}$ when $i + 1$ is even.
  It is clear that the solution determined by the symplectic condition
  is meromorphic, and hence by the above analysis, it is also
  holomorphic.
\end{proof}

Let $R_1$ and $R_2$ be as in the lemma.
We can then define a convergent CohFT $\Omega^1$ defined outside the
discriminant locus by the $R$-matrix action of $R_1$ on the trivial
CohFT.
We will show that $\Omega^1$ extends to the discriminant locus by
comparing it to the (everywhere defined) convergent CohFT $\Omega^3$
defined by the $R$-matrix action $\Omega^3 = (R_1 R_2^{-1}) \Omega^2$.
The CohFTs $\Omega^3$ and $\Omega^1$ are very similar but the
underlying trivial theories do not agree, in particular the units
$\mathbf 1_{\Omega^1}$ and $\mathbf 1_{\Omega^2}$ as well as the norms
$\Delta_{1i}^{-1}$ and $\Delta_{2i}^{-1}$ of the idempotents are
different.
We will fix this using a shift.

For this, recall the description \eqref{eq:reconstrb} of the
reconstruction using the basis of normalized idempotents.
A local contribution at a vertex of color $i$ for the reconstruction
of $\Omega^3$ is of the form
\begin{equation*}
  \Delta_{2i}^{\frac{2g - 2 + n}2} \sum_{k = 0}^\infty \frac{\Delta_{2i}^{\frac k2}}{k!} \pi_*\left(\prod_{j = 1}^n \alpha_j^i \prod_{j = 1}^k \psi_j (\Id^i - (R_1^{-1}(\psi_j))^i)\mathbf 1_{\Omega^2}\right),
\end{equation*}
where $\pi$ forgets the last $k$ markings and $\alpha_j^i$ are some
formal series in $\psi_j$ whose coefficients are holomorphic functions
on the Frobenius manifold.
To circumvent convergence issues, let $v$ be a formal flat vector
field and let $v_\mu$ and $v_i$ be the coordinates of $v$ when written
in a basis of flat coordinates or in the basis of normalized
idempotents, respectively.
We can further modify the CohFT by shifting along $v\psi$:
\begin{equation}
  \label{eq:dilshift}
  \Omega_{g, n}^4 (\alpha_1, \dotsc, \alpha_n)
  := \sum_{k = 0}^\infty \frac 1{k!} \pi_* \Omega_{g, n + k}^3 (\alpha_1, \dotsc, \alpha_n, \psi v, \dotsc, \psi v)
\end{equation}
We obtain a well-defined convergent CohFT defined over the ring of
power series in the $v_\mu$.
For $\Omega^4$ the local contribution at a vertex of color $i$ is
\begin{equation*}
  \Delta_{2i}^{\frac{2g - 2 + n}2} \sum_{k = 0}^\infty \frac{\Delta_{2i}^{\frac k2}}{k!} \pi_*\left(\prod_{j = 1}^n \alpha_j^i \prod_{j = 1}^k \psi_j \left[\Id^i\mathbf 1_{\Omega^2} - (R_1^{-1}(\psi_j))^i(\mathbf 1_{\Omega^2} - v)\right]\right).
\end{equation*}

Recall that the dilaton equation implies that
\begin{equation*}
  \frac 1{(1 - a)^{2g - 2 + n}}
  = \sum_{k = 0}^\infty \frac 1{k!} \pi_* \left(\prod_{j = 1}^k a\psi_j\right),
\end{equation*}
where $a$ is a formal variable, or equivalently
\begin{equation*}
  1 = \sum_{k = 0}^\infty \frac{(1 + b)^{2 - 2g - n - k}}{k!} \pi_* \left(\prod_{j = 1}^k b\psi_j\right),
\end{equation*}
where $b=a/(1-a)$.
We will apply this identity locally at every vertex.
At a vertex of color $i$ we use $-\sqrt{\Delta_{2i}}v_i$ for $b$.
Then the local contribution at a vertex of color $i$ is
\begin{equation*}
  \Delta_{3i}^{\frac{2g - 2 + n}2} \sum_{k = 0}^\infty \frac{\Delta_{3i}^{\frac k2}}{k!} \pi_*\left(\prod_{j = 1}^n \alpha_j^i \prod_{j = 1}^k \psi_j \left[\Id^i - (R_1^{-1}(\psi_j))^i\right](\mathbf 1_{\Omega^2} - v)\right),
\end{equation*}
where
\begin{equation*}
  \Delta_{3i}^{-1/2} = \Delta_{2i}^{-1/2} - v_i.
\end{equation*}
Notice that now (again) the sum in $k$ is finite in each cohomological
degree.
Therefore we can specialize $v$.
We will take $v$ to be the vector
$\mathbf 1_{\Omega^2} - \mathbf 1_{\Omega^1}$, which is holomorphic by
Lemma~\ref{lem:basischange}, and thus
$v_i = \Delta_{2i}^{-1/2} - \Delta_{1i}^{-1/2}$.
In this case, the $\Delta_{3i}$ specialize to $\Delta_{1i}$.
We have therefore arrived exactly at the reconstruction formula for
$\Omega^1$.
In particular, with the specialization of $v$, $\Omega^4$ is the same
as $\Omega^1$ and therefore $\Omega^1$ is also holomorphic along the
discriminant.
Hence $\Omega^1$ is a suitable local extension of the Frobenius
manifold we started with to a convergent CohFT.
We note that the extension in not unique but depends on the choice of
integration constants in Lemma~\ref{lem:Rreg}.

\subsection{Equivalence of relations}
\label{sec:cohft:eq}

We will prove Theorem~\ref{thm:main} in this section.

First notice that the dimension of a convergent CohFT $\Omega$ can be
increased by one without changing the set of relations by the
construction of Example~\ref{ex:extend}.
So we can assume that the CohFTs we are trying to compare have the
same dimension.

Next, recall from Section~\ref{sec:cohft:taut} that the tautological
relations of a semisimple convergent CohFT $\Omega$ are defined via
coefficients of the part of the Givental--Teleman classification
singular in the discriminant.
Therefore the relations do not change when removing the codimension
two set of singular points of the discriminant from the Frobenius
manifold underlying $\Omega$.

In order to prove Theorem~\ref{thm:main}, it is therefore enough to
show that the relations coincide for two semisimple, equal dimensional
convergent CohFTs $\Omega^1$, $\Omega^2$ such that each Frobenius
manifold contains a smooth point of the discriminant and is small
enough for Theorem~\ref{thm:localext} to apply directly to $\Omega^1$.

By the proof of Theorem~\ref{thm:localext}, an extension of the
Frobenius manifold underlying $\Omega^1$ to a CohFT can be constructed
from $\Omega^2$ by a holomorphic $R$-matrix and a holomorphic shift.
To prove Theorem~\ref{thm:main}, it therefore suffices to show that
these two operations preserve tautological relations and that the
integration constants in the proof of Lemma~\ref{lem:Rreg} can
be chosen such that the constructed CohFT coincides with $\Omega^1$.
We now prove these statements.
\begin{lemma}
  \label{lem:preserve}
  The $R$-matrix action by a holomorphic $R$-matrix preserves tautological
  relations.
\end{lemma}
\begin{proof}
  Let $\Omega'$ be obtained from $\Omega$ from the $R$-matrix action
  of $R$.
  Then in the description of the $R$-matrix action in
  Section~\ref{sec:cohft:reconstr}, all arguments are holomorphic
  vector fields on the Frobenius manifold with values in power series
  in $\psi$-classes.
  $\Omega'_{g, n}$ in each cohomological degree is obtained by a
  finite sum of push-forwards under the gluing map of products of
  $\Omega$ (with possibly additional markings) multiplied by monomials
  in $\psi$ classes and with holomorphic vector fields as arguments.
  Therefore any singularities of the reconstruction of $\Omega'$ are
  the result of singularities in the reconstruction of $\Omega$.
  So we can write the relations $q(\widetilde{\Omega'}_{g, n})$ in
  terms of the relations from $\Omega$ as of
  Definition~\ref{def:rels}.
  Since the ideal of relations of $\Omega'$ is stable under the
  tautological maps, we can also express any tautological relation of
  $\Omega'$ as in Definition~\ref{def:rels} in terms of relations from
  $\Omega$.

  Since $R$-matrices are power series starting with the identity
  matrix, by using $R^{-1}$, we can also write the relations of
  $\Omega$ in terms of relations from $\Omega'$.
\end{proof}

The shift-construction \eqref{eq:dilshift} clearly expresses any
relation from $\Omega^4_{g, n}$ in terms of relations from
$\Omega^3_{g, n + m}$ for various $m \ge 0$.

Now, we finally show that taking the $R$-matrix of $\Omega^1$ is a
suitable choice for $R_1$ in Lemma~\ref{lem:Rreg}.
We will argue that otherwise $\Omega^1$ or $\Omega^2$ would not be
defined at the discriminant.

For simplicity, we will make use of the following stability result.
It should also be possible to use estimates or congruence properties
of intersection numbers instead.
\begin{theorem}[Boldsen \cite{Bo12}, Looijenga \cite{Lo96}]
  \label{thm:bolo}
  Let $M_{g, n} \subset \Mbar_{g, n}$ be the moduli space of smooth
  curves.
  Then, for any $k < \frac g3$ the vector space $H^{2k}(M_{g, n})$ is
  freely generated by the set of monomials of cohomological degree
  $2k$ in the classes
  $\kappa_1, \dotsc, \kappa_k, \psi_1, \dotsc, \psi_n$.
\end{theorem}

We use the local coordinates $t, t_0, u_{\ge 3}$ from the proof of
Lemma~\ref{lem:Rreg}.
Let $i$ be the lowest degree in $z$ where $R_1 R_2^{-1}$ is not
holomorphic.
By the proof of Lemma~\ref{lem:Rreg}, the non-holomorphic part in this
degree is a constant multiple of the block-diagonal matrix with
upper-left block
\begin{equation*}
  \begin{pmatrix}
    0 & t^{1/2} \\
    t^{-1/2} & 0
  \end{pmatrix}
\end{equation*}
and zeros everywhere else.
By Theorem~\ref{thm:bolo}, we can consider the $\psi_1^i$-coefficient
of
\begin{equation*}
  \Omega_{g, 1}^1 \left(\frac\partial{\partial t_0}\right)\Big|_{M_{g, 1}} - \Omega_{g, 1}^2 \left(\frac\partial{\partial t_0}\right)\Big|_{M_{g, 1}}
\end{equation*}
for large $g$.
Its lowest order term when expanded in $t$ is up to a nonzero factor
given by
\begin{multline*}
  t^{-\frac 12} \sqrt t \left(\sqrt{2\sqrt t}^{2g - 2 + 1 - 1} - \sqrt{-2\sqrt t}^{2g - 2 + 1 - 1}\right) \\
  = 2^{g - 1} \left((\sqrt t)^{g - 1} - (-\sqrt t)^{g - 1}\right),
\end{multline*}
and therefore not holomorphic in $t$ for even $g$.
This is a contradiction to the assumption that both $\Omega^1$ and
$\Omega^2$ are defined on the discriminant locus.

\subsection{Global extension}
\label{sec:cohft:glocalext}

The local extension Theorem~\ref{thm:localext} leaves open the
question under what conditions a semisimple Frobenius manifold can
(globally) be extended to a CohFT.
In Section~\ref{sec:ex:global}, we will see that the restrictions put
on integration constants of the $R$-matrices in Lemma~\ref{lem:Rreg}
do not always fit together globally.
\begin{conjecture}
  \label{conj:globalext}
  Let $M$ be an $N$-dimensional semisimple Frobenius manifold such
  that it possesses a holomorphic genus one potential $\mathrm dG$.
  Then there exists a convergent CohFT with underlying Frobenius
  manifold $M$.
\end{conjecture}

On the other hand, when the Frobenius manifold is conformal, such an
extension to a CohFT exists by the following simple argument.
There is a unique homogeneous solution to the flatness equation
\eqref{eq:RDE} and by construction it is meromorphic along the
discriminant.
Since by Lemma~\ref{lem:Rreg} all possible solutions are either
holomorphic in the discriminant or are multivalued, the homogeneous
solution has in fact to be holomorphic.

\subsection{Examples}
\label{sec:cohft:ex}

\subsubsection{Extending the comparison to non-smooth points on the
  discriminant}

We want to illustrate how the comparison between relations in the
proof of Theorem~\ref{thm:main} via an identification of coordinates
and vector fields, an $R$-matrix action and a shift, does not directly
extend to give a way to explicitly write the relations near a singular
point of the discriminant in terms of the $A_2$- (3-spin) relations.

Let us consider the comparison between the $A_2 \times A_1$ and $A_3$
singularities.
We will see that already the identification between points and vector
fields behaves badly.
This is the simplest example we can consider since in two dimensions
the discriminant locus is a union of parallel lines and in particular
is non-singular.

The Frobenius manifold of the $A_3$-singularity $x^4/4 = 0$ is based on
the versal deformation space
\begin{equation*}
  f(x) = \frac {x^4}4 + t_2 x^2 + t_1 x + t_0.
\end{equation*}
Here $t_0$, $t_1$ and $t_2$ are coordinates on the Frobenius manifold.
The ring structure is given by the Milnor ring
\begin{equation*}
  \CC[t_0, t_1, t_2][x] / f'(x),
\end{equation*}
where $x = \frac\partial{\partial t_1}$.
The discriminant of the minimal polynomial $f'$ of $x$ is
$-32 t_2^3 - 27 t_1^2$ and therefore the discriminant locus has a cusp
at $t_1 = t_2 = 0$.
The metric is in the basis $\{1, x, x^2\}$ given by
\begin{equation*}
  \begin{pmatrix}
    0 & 0 & 1 \\
    0 & 1 & 0 \\
    1 & 0 & -2t_2
  \end{pmatrix}.
\end{equation*}
Therefore the basis $\{1, x, x^2\}$ is flat up to a determinant one
basis change.

We go to a sixfold ramified cover of the Frobenius manifold on which
we can define the critical points $\zeta_1, \zeta_2, \zeta_3$ of
$f(x)$ as holomorphic functions.
Let $u_1, u_2, u_3$ be the corresponding critical values.
Part of the discriminant locus is described by the equation
$\zeta_1 = \zeta_2$.
Locally we use $\phi := \zeta_1 - \zeta_2$, $\zeta_3$ and $t_0$ as new
coordinates.
Reexpressing in terms of these coordinates gives
\begin{align*}
  \zeta_1 &= -\frac 12 \zeta_3 + \frac 12 \phi, \\
  \zeta_2 &= -\frac 12 \zeta_3 - \frac 12 \phi, \\
  u_1 &= t_0 + \frac 3{64} \zeta_3^4 - \frac 5{32} \zeta_3^2 \phi^2 + \frac 18 \zeta_3 \phi^3 - \frac 1{64} \phi^4, \\
  u_2 &= t_0 + \frac 3{64} \zeta_3^4 - \frac 5{32} \zeta_3^2 \phi^2 - \frac 18 \zeta_3 \phi^3 - \frac 1{64} \phi^4, \\
  u_1 - u_2 &= \frac 14 \zeta_3 \phi^3, \\
  u_3 &= t_0 - \frac 38 \zeta_3^4 + \frac 18 \zeta_3^2 \phi^2.
\end{align*}
The idempotents are given by
\begin{align*}
  \frac\partial{\partial u_1} =& \frac{(x - \zeta_2)(x - \zeta_3)}{(\zeta_1 - \zeta_2)(\zeta_1 - \zeta_3)}
  = \frac{x^2 + (-\frac 12 \zeta_3 + \frac 12 \phi)x - \frac 12 \zeta_3^2 - \frac 12 \zeta_3 \phi}{-\frac 32 \zeta_3 \phi+ \frac 12 \phi^2}, \\
  \frac\partial{\partial u_2} =& \frac{(x - \zeta_1)(x - \zeta_3)}{(\zeta_2 - \zeta_1)(\zeta_2 - \zeta_3)}
  = \frac{x^2 + (-\frac 12 \zeta_3 - \frac 12 \phi)x - \frac 12 \zeta_3^2 + \frac 12 \zeta_3 \phi}{\frac 32 \zeta_3 \phi+ \frac 12 \phi^2}, \\
  \frac\partial{\partial u_3} =& \frac{(x - \zeta_1)(x - \zeta_2)}{(\zeta_3 - \zeta_1)(\zeta_3 - \zeta_2)}
  = \frac{x^2 + \zeta_3 x + \frac 14 \zeta_3^2 - \frac 14 \phi^2}{\frac 94 \zeta_3^2 - \frac 14 \phi^2}
\end{align*}
so that they become after normalization
\begin{multline*}
  \frac{x^2 + (-\frac 12 \zeta_3 + \frac 12 \phi)x - \frac 12
    \zeta_3^2 - \frac 12 \zeta_3 \phi}{\sqrt{-\frac 32 \zeta_3 \phi+
      \frac 12 \phi^2}},
  \frac{x^2 + (-\frac 12 \zeta_3 - \frac 12 \phi)x - \frac 12 \zeta_3^2 + \frac 12 \zeta_3 \phi}{\sqrt{\frac 32 \zeta_3 \phi+ \frac 12 \phi^2}}, \\
  \frac{x^2 + \zeta_3 x + \frac 14 \zeta_3^2 - \frac 14
    \phi^2}{\sqrt{\frac 94 \zeta_3^2 - \frac 14 \phi^2}}.
\end{multline*}

For $A_2 \times A_1$, let us assume that $u_3$ corresponds to the
$A_1$-direction and that the norm of the idempotent in that direction
is one.
We can write
\begin{align*}
  u_1 = x_0 - \frac 23 (-x_1)^{3/2}, \quad u_2 = x_0 + \frac 23 (-x_1)^{3/2},
\end{align*}
where $x_0$ and $x_1$ are flat coordinates corresponding to $t_0$ and
$t_1$ in Example~\ref{ex:3spin}.

We should therefore identify
\begin{equation*}
  \phi \stackrel{!}{=} -2\left(\frac 23\right)^{1/3} \zeta_3^{-1/3} \sqrt{-x_1}.
\end{equation*}

Let us consider how we identify the $A_3$-singularity basis
$\{1, x, x^2\}$ and the flat basis of $A_2 \times A_1$ via the
identification of their normalized idempotents.
If we write the identity of $A_2 \times A_1$ in terms of
$\{1, x, x^2\}$, the $x^2$-coefficient is
\begin{multline*}
  \frac{\sqrt{2 \sqrt{-x_1}}}{\sqrt{-\frac 32 \zeta_3 \phi+
      \frac 12 \phi^2}} + \frac{\sqrt{-2 \sqrt{-x_1}}}{\sqrt{\frac 32 \zeta_3 \phi+
      \frac 12 \phi^2}} + \frac 1{\sqrt{\frac 94 \zeta_3^2 - \frac 14 \phi^2}} \\
  = \frac{2\zeta_3^{1/6}}{\sqrt{-3c \zeta_3 + c \phi}} + \frac{2\zeta_3^{1/6}}{\sqrt{3c \zeta_3 + c \phi}} + \frac 2{\sqrt{9\zeta_3^2 - \phi^2}},
\end{multline*}
where
\begin{equation*}
  c = -2\left(\frac 23\right)^{1/3}.
\end{equation*}
The coefficient is well-defined on generic points of the part
$\zeta_1 = \zeta_2$ ($\phi = 0$) of the discriminant, but when fixing
some $\phi \neq 0$ the function has a singularity at $\zeta_3 = 0$.

\subsubsection{Obstructions to extending $R$-matrices}
\label{sec:ex:global}

We want consider the class of two-dimensional Frobenius manifolds
with flat coordinates $t_0, t$, flat metric
\begin{equation*}
  \eta =
  \begin{pmatrix}
    0 & 1 \\
    1 & 0
  \end{pmatrix}
\end{equation*}
and quantum product
\begin{equation*}
  \left(\frac\partial{\partial t}\right)^2 = f \frac\partial{\partial t_0}
\end{equation*}
for a holomorphic function $f(t)$.
The corresponding Gromov--Witten potential is
\begin{equation*}
  \frac 12 t_0^2 t + F,
\end{equation*}
where $F(t)$ is a third anti-derivative of $f(t)$.

The differential equation satisfied by the $R$-matrix in flat
coordinates can be made explicit:
\begin{equation}
  \label{eq:Rex}
  \left[R,
    \begin{pmatrix}
      0 & f \\
      1 & 0
    \end{pmatrix}
  \right] + z \dot R + z \frac{\dot f}{4f}
  \begin{pmatrix}
    -1 & 0 \\
    0 & 1
  \end{pmatrix}
  R = 0
\end{equation}
We first want to show that for any solution $R$, the $z^1$-coefficient
is not holomorphic for all $f$.
For this we set
\begin{equation*}
  R =
  \begin{pmatrix}
    1 + az & 0 + bz \\
    0 + cz & 1 + dz
  \end{pmatrix}
  + O(z^2).
\end{equation*}
From \eqref{eq:Rex} in degree $z^1$ we obtain
\begin{equation*}
  b - fc - \frac{\dot f}{4f} = 0, \qquad a = d.
\end{equation*}
From \eqref{eq:Rex} in degree $z^2$, we see that $a = d$ is an
integration constant.
We obtain an interesting differential equation for $c$:
\begin{equation*}
  2f \dot c + \dot f c + \frac{\ddot f}{4f} - \frac{5\dot f^2}{16f^2} = 0
\end{equation*}
If we substitute
\begin{equation*}
  c = \frac{\gamma}f - \frac 5{48} \frac{\dot f}{f^2},
\end{equation*}
it becomes
\begin{equation*}
  2\dot\gamma - \frac{\dot f}f \gamma + \frac{\ddot f}{24f} = 0.
\end{equation*}
So $\gamma$ is determined up to a multiple of a root of $f$ and in
particular, if $f$ has somewhere a simple zero, there exists at most
one solution meromorphic on all of $\CC^2$.

If $f$ is linear, $\gamma = 0$ is clearly a holomorphic solution.
If $f$ is quadratic with non-vanishing discriminant, there is still a
holomorphic solution.
For example for
\begin{equation*}
  f(t) = t(t + 1)
\end{equation*}
the solution is
\begin{equation*}
  \gamma = \frac t6 + \frac 1{12}.
\end{equation*}
In larger degree, we stop having meromorphic solutions.
In the example
\begin{equation*}
  f(t) = t(t^2 - 1),
\end{equation*}
after substituting
\begin{equation*}
  \gamma = f \delta + \frac t8,
\end{equation*}
we arrive at the differential equation
\begin{equation*}
  (t^2 - 1) 2t\dot \delta + (3t^2 - 1)\delta + \frac 18 = 0.
\end{equation*}
We see that $\delta$ is meromorphic in $t$ if and only it is so in
$u := t^2$.
In the new variable the differential equation is
\begin{equation*}
  4u(u - 1) \delta' + (3u - 1)\delta + \frac 18 = 0.
\end{equation*}
From generic semisimplicity we also know that $\delta$ has to be
holomorphic except for $u = 0$ and $u = 1$.
Around $u = 0$ and $u = 1$ there are unique meromorphic solutions
\begin{equation*}
  \frac 18 \sum_{i = 0}^\infty \frac{4i + 3}{4i + 1} u^i, \qquad
  -\frac 1{16} \sum_{i = 0}^\infty \frac{4i + 3}{4i + 2} (1 - u)^i,
\end{equation*}
but these obviously do not agree (consider $u = \frac 12$).

We now want to check that the corresponding genus one potential will
also be singular, in agreement with Conjecture~\ref{conj:globalext}.
For this, we look at the case when $a = d = 0$, and compute the
codimension one part of the reconstructed CohFT on $\Mbar_{1, 1}$ with
an $\frac\partial{\partial t}$-insertion.
From the trivial graph, we obtain the contribution
\begin{equation*}
  -2\left(\gamma + \frac 7{48} \frac{\dot f}f\right) \psi_1 + 2\left(\gamma - \frac 5{48} \frac{\dot f}f\right) \kappa_1
\end{equation*}
and in addition we have the contribution
\begin{equation*}
  2\gamma + \frac 2{48} \frac{\dot f}f
\end{equation*}
of the irreducible divisor $\delta_0$.
From
\begin{equation*}
  \int_{\Mbar_{1, 1}} \psi_1 = \int_{\Mbar_{1, 1}} \kappa_1 = \frac 1{12} \int_{\Mbar_{1, 1}} \delta_0 = \frac 1{24}
\end{equation*}
we see that the correlator equals $\gamma$, which is not holomorphic
on all of the Frobenius manifold.

\subsection{Other relations from cohomological field
  theories}
\label{sec:cohft:gw}

For a convergent CohFT depending on additional parameters, there are
ways different from Definition~\ref{def:rels} to obtain tautological
relations from the reconstruction of semisimple CohFTs.
We will consider here the example of the equivariant Gromov--Witten
theory of a toric variety, which is dependent on equivariant and
Novikov parameters.
We will see that the additional relations are still a consequence of
Pixton's relations.

Let $T = (\CC^*)^m$ and let
$H_T^*(\mathrm{pt}) = H^*(BT) = \CC[\lambda_1, \dotsc, \lambda_m]$ be
the $T$-equivariant cohomology ring of a point.
Let $X$ be an $m$-dimensional smooth, toric variety with a basis
$\{p_1, \dotsc, p_N\}$ of its cohomology, which we can also lift to a
basis in $T$-equivariant cohomology.
Let $\beta_1, \dotsc, \beta_N$ be the dual homology basis.
The Novikov ring is a completion of the semigroup ring of effective
classes $\beta \in H_2(X; \ZZ)$.
We use $q^\beta$ to denote the generator corresponding to a
$\beta \in H_2(X; \ZZ)$.

A family of $N$-dimensional CohFTs on the state space $H_{\CC^*}^*(X)$
can be defined by setting
\begin{equation*}
  \Omega_{g, n}(\alpha_1, \dotsc, \alpha_n) = \sum_\beta q^\beta p_*\left(\prod_{i = 1}^n \ev_i^*(\alpha_i) \cap [\Mbar_{g, n}(X; \beta)]^{vir}\right),
\end{equation*}
where the sum ranges over all effective classes
$\beta \in H_2(X; \ZZ)$, $p$ is the projection from the moduli space
of stable maps to $\Mbar_{g, n}$, and $\ev_i$ is the $i$th evaluation
map.
At this point $\Omega_{g, n}$ is just a power series in the Novikov
variables.

Using virtual localization with respect to the $T$-action on $X$, the
CohFT $\Omega_{g, n}$ can be effectively computed.
The corresponding $R$-matrix is however only defined after localizing
the equivariant parameters $\lambda_1, \dotsc, \lambda_m$.
This gives rise to tautological relations defined similarly to those
of Definition~\ref{def:rels}.
The reconstruction formula (which in this case is a consequence of
virtual localization) gives an expression
$\widetilde\Omega^{\mathrm{loc}}_{g, n}$ for $\Omega_{g, n}$ which
does not seem to admit the non-equivariant limit $\lambda_i \to 0$.
On the other hand, the well-defined non-equivariant Gromov--Witten
theory of $X$ must be recovered in this limit.
The necessary cancellation of poles in the equivariant parameters
gives rise to tautological relations.
In the example of $X = \PP^1$, the resultant relations actually imply
Pixton's relations \cite{Ja17, Ja13P}.

We now compare these relations from equivariant localization to the
relations considered in Definition~\ref{def:rels}.
From results \cite{Ir07} of Iritani, it follows that the sum over
$\beta$ defining $\Omega_{g, n}$ converges in a neighborhood of the
origin, and that the CohFT induces a convergent CohFT varying
analytically on the additional Novikov and equivariant parameters.
We can therefore view $\Omega_{g, n}$ as valued on the space of
functions $\cO_U$ on an open neighborhood $U$ of the orgin in the
product of the Frobenius manifold and the space of possible Novikov
and equivariant parameters.
By mirror symmetry (see \cite[Section~6]{Ir07}), it follows also that
the coefficients of the $R$-matrix are valued in $\cO_{U \setminus D}$
where $D \subset U$ is a discriminant locus.
Therefore the reconstruction gives elements
$\widetilde\Omega_{g, n} \in \mathcal S_{g, n} \otimes (V^*)^{\otimes
  n} \otimes \cO_{U \setminus D}$.
As in Definition~\ref{def:rels}, projecting to
$\cO_{U \setminus D}/\cO_U$ yields tautological relations.
Since for any value of the Novikov and equivariant parameters, the
CohFT is generically semisimple, we can apply Theorem~\ref{thm:main}
to see that all these tautological relations are consequences of
Pixton's relations.
Restricting $\widetilde\Omega_{g, n}$ to the origin on the Frobenius
manifold and Taylor expanding in the Novikov parameters, we recover
$\widetilde\Omega^{\mathrm{loc}}_{g, n}$.
It follows that the tautological relations from considering the
non-equivariant limit are also consequences of Pixton's relations.

\begin{remark}
  A similar strategy should also work for toric orbifolds.
  The special case of $\PP^1$ with two orbifold points is used in
  \cite{ClJa16P}.
\end{remark}

\bibliographystyle{plain}
\bibliography{cohftrel}
\addcontentsline{toc}{section}{References}

\vspace{+8 pt}
\noindent
Department of Mathematics \\
University of Michigan \\
2074 East Hall \\
530 Church Street \\
Ann Arbor \\
MI-48109 \\
USA

\end{document}